\documentclass[11pt,leqno]{article}

\usepackage{amsthm}
\usepackage{amsmath}
\usepackage{amssymb}
\usepackage{amsfonts}
\usepackage{amscd}
\usepackage{palatino}
\usepackage[mathscr]{eucal}
\usepackage{epsfig}
\usepackage{verbatim}
\usepackage{latexsym}
\usepackage{graphics}

\begin{document}

\newcounter{chn}
\newenvironment{nch}[1]
               {\addtocounter{chn}{1}
               \newpage
                      \begin{center}
                      {\Large {\bf {#1}}}
                      \end{center}
                      \vspace{5mm}}{}

%The Theorem env.
\newcounter{thn}
\renewenvironment{th}
{\addtocounter{thn}{1}
 \em 
     \noindent 
     {\sc Theorem \arabic{thn} . ---}}{}

%Theorem without numbers
\newenvironment{theo}
   {\em
        \vspace{2mm}
        \noindent
                {\sc Theorem\/. ---}}{}

%Lemma without numbers
\newenvironment{lem}
   {\em
        \vspace{2mm}
        \noindent
                {\sc Lemma\/. ---}}{}

%Prop without numbers
\newenvironment{pro}
   {\em
        \vspace{2mm}
        \noindent
                {\sc Proposition\/. ---}}{}

%The Lemma env.
\newcounter{len}
\renewenvironment{le}
{\addtocounter{len}{1}
 \em 
     \vspace{2mm}
     \noindent 
     {\sc Lemma \arabic{len}. ---}}{}

%The Proposition env.
\newcounter{prn}
\newenvironment{prop}
{\addtocounter{prn}{1}
  \em 
      \vspace{2mm}
      \noindent 
      {\sc Proposition \arabic{prn}. ---}}{}

\newcounter{exn}
\newenvironment{prob}
                     {\addtocounter{exn}{1} 
                      \vspace{2mm}
                      \noindent 
                      {\em Exercise \arabic{exn} :}}{}

\newcounter{corn}
\newenvironment{cor}
  {\em 
     \addtocounter{corn}{1}
     \vspace{2mm}   
     \noindent 
        {\sc Corollary \arabic{corn}. ---}}{\vspace{1mm}}

%Corollary without numbers
\newenvironment{cor1}
   {\em 
        \vspace{2mm}
        \noindent
                {\sc Corollary\/. ---}}{}

%Definitions
\newcounter{defn}
\newenvironment{df}{
                    \vspace{2mm}
                    \addtocounter{defn}{1} 
                    \noindent 
                            {\sc Definition \arabic{defn}. ---}}{\vspace{1mm} }{}
%paragraphs
\newcounter{pgn}
\newenvironment{pg}{
                    \vspace{2mm}
                    \addtocounter{pgn}{1}
                    \noindent
                            {\arabic{scn}.\arabic{pgn}\;}}{}

%The section env.
\newcounter{scn}
\newenvironment{nsc}[1]
              {\addtocounter{scn}{1}
              \vspace{3mm}
                     \begin{center}
                     {\sc { \arabic{scn}. {#1}}}
                     \end{center}
              \setcounter{equation}{0}
              \setcounter{prn}{0}
              \setcounter{len}{0}
              \setcounter{thn}{0}
              \setcounter{pgn}{0}
              \setcounter{corn}{0}
              \vspace{3mm}}{}

\begin{center}

\begin{center}
\Large {\bf 

{

The Large Sieve Inequality for Quadratic Polynomial Amplitudes

}
}

\end{center}
\vspace{5mm}
\end{center}

\vspace{-10mm}
\begin{center}

\it{

by \\

\vspace{3mm}

Gyan Prakash and D.S. Ramana 

}
\end{center}

\begin{nsc}{Introduction}

\vspace{2mm}
\noindent
An important requirement in the context of inequalities of the large sieve type  is to obtain estimates for the sum $\sum_{x \in \mathcal{X}}  |\sum_{i \in I } a_i e(xy_i))|^2$, where $\mathcal{X}$ is a well-spaced set of real numbers, $I$ is a finite set, $\{a_i\}_{i \in I}$ are complex numbers 
and $ \{y_i\}_{i\in I}$ is a {\em sparse subsequence} of the integers. 

\vspace{2mm}
\noindent
Basic examples of sparse sequences of integers are provided by the sequence of
values of polynomials of degree $\geq 2$ with integer coefficients. The
present article is concerned with the case when polynomial is of degree
2. Indeed, in a recent note, Liangyi Zhao \cite{zh}, showed, by an elegant application of the double large sieve
inequality of Bombieri and Iwaniec, that one has the estimate given below,
which we state with the aid of the following notation.

\vspace{2mm}
\noindent
When $Q$ is a real number $\geq 1$, the Farey series of order $Q$ is the sequence of rational numbers in $(0,1)$ whose denominators, when expressed lowest form, do not exceed $Q$. Given a finite sequence of complex numbers $\{a_i\}_{i \in I}$, indexed by a finite set $I$, we write $\|a\|^2$ to denote $\sum_{i \in I} |a_i|^2$.

\vspace{2mm}
\noindent
{\sc Theorem (L. Zhao).--- }{\em Let $Q$ be a real number $\geq 1$ and suppose that ${\mathcal F}(Q)$ is the Farey sequence of order $Q$ and that $P(T) = c_0T^2 + c_1 T + c_2$ is a polynomial in ${\bf R}[T]$ with $c_0 \neq 0$, $c_1/c_0 = p/q \in {\bf Q} $, $c_1 > 0$ and $(p,q) =1$. When $\{a_i\}_{M < i \leq M+N}$ is a sequence of complex numbers indexed by the integers in the interval $(M, M+N]$, where $M$ and $N$ are integers with $N \geq 1$, we have 

\begin{equation}
\label{lsth}
\sum_{x \in \mathcal{{\mathcal F}(Q)}}  \left| \sum_{M < i \leq M+N } a_i e( x P(i))\right|^2  \; \ll \; (Q^2 + Q \sqrt{c_0 N (|M| + 2N + p/q + 1)} \Pi \|a\|^2 \; ,
\end{equation}

\noindent
where the implied constant depends on $\epsilon$ alone and 

\[
\Pi = \left(\frac{q}{c_0} + 1\right)^{\frac{1}{2} + \epsilon} [N q(|M| + N) + |p| + q/c_0]^{\epsilon} \; .
\]

}

\vspace{2mm}
\noindent
Zhao has shown in \cite{zh}, Section 3, page 4 that one may device examples of $P(T)$  and the sequence $\{a_i\}$ for which the left hand side of~\eqref{lsth}
is $\gg (NQ + Q^2)\|a\|^2$, thereby showing that \eqref{lsth} is essentially the best possible.

\vspace{2mm}
\noindent
There, however, remains the question of determining the extent to which the dependence of the right hand side of \eqref{lsth} on $M$ and the coefficients of $P(T)$ may be reduced, given that the trivial estimate for the left hand side of \eqref{lsth} is independent of these parameters.
In the present article we combine Zhao's method
 with an interpolation argument due to Heath-Brown \cite{Hb}
 to show that atleast a modest improvement upon the above theorem along these lines is certainly possible. More precisely,  the corollary to Theorem 1 of Section 3 below shows that 
, under the same hypotheses as in the above theorem, (\ref{lsth}) holds with the factor $\Pi$ replaced  with

\[
\Pi^{\prime} = \pi \left(\frac{2q}{c_0} + 1\right)^{\frac{1}{2}} \sup_{1 \leq n \leq 144N^4} r(n) \; ,   
\]                                              

\noindent
and the $\ll$ in (\ref{lsth}) replaced by $\leq$. Here $r(n)$ , for an integer
$n$, is the number of integer points on $(x,y)$ satisfying $x^2 + y^2 = n$. 

\vspace{2mm}
\noindent
The reader will note with interest that S. Baier has already shown in \cite{st} that an analog of our conclusion recorded in the corollary to Theorem 1 of Section 3 holds for {\em all} quadratic polynomials with {\em real coefficients} when one replaces $\|a\|^{2}$ with ${\rm Card}(I) \sup_{i \in I} |a_i|^2$. Moerover, this work of Baier also investigates what might be expected to hold for higher degree polynomials in this context.

\end{nsc}

\begin{nsc}{Preliminaries}
\vspace{2mm}
\noindent
{\em 2.1. $\delta$-Spaced Sets. ---}  Let $\delta > 0$. A $\delta$-spaced set of real numbers is a  finite set $\mathcal{X}$ of  distinct real numbers containing atleast two elements  and such that $|x-x^{\prime}| \geq \delta$, whenever $x$ and $x^{\prime}$ are distinct elements of $\mathcal{X}$. Let $\mathcal{X}$ be a $\delta$- spaced set and $\epsilon > 0$. We then set $S(\epsilon, x) = {\rm{Card}}\left( x^{\prime} \in {\mathcal{X}}\, | \; |x - x^{\prime}| \leq \epsilon \right)$. We have

\vspace{-3mm}
\begin{equation}
\label{spacing}
S(\epsilon, x) \; \leq \;
1 + \frac{2\epsilon}{\delta} \;\; \text{for all $x$ in $\mathcal{X}$.}
\end{equation}

\vspace{3mm}
\noindent
{\em 2.2. The function ${\phi}$. ---} Let $\phi(t)$ denote the characteristic function of the interval $[-1/2, 1/2]$ and for each $\epsilon > 0$, set $\phi_{\epsilon}(t) = \phi(t/2\epsilon)$. Thus $\phi_{\epsilon}(t)$ is the characteristic function of the interval $[-\epsilon, \epsilon]$. On setting $\sin{t}/t$ to  1 when $t =0$ we have 

\vspace{-3mm}
\begin{equation}
\label{phi}
\widehat{\phi_{\epsilon}}(t) = 2 \epsilon \widehat{\phi}(2 \epsilon t) = 2\epsilon \, \left( \frac{\sin{2 \pi \epsilon t}}{2 \pi \epsilon t} \right) \; .
\end{equation}

\noindent
The classical inequality $2/\pi \leq \sin t/ t \leq 1$ for $t$ in $[0, \pi/2]$ then implies the inequality 

\vspace{-3mm}
\begin{equation}
\label{phi1}
\frac{1}{2\epsilon}\; \leq \; \frac{1}{\widehat{\phi_{\epsilon}}(t)} \; \leq \; \frac{\pi}{4 \epsilon} \; ,
\end{equation}

\noindent
for all $t$ with $\epsilon t$ in $[-1/4, 1/4]$.

\vspace{2mm}
\noindent
{\em 2.3. A Simple Majorisation Principle. ---} Let $I$ be a finite set, $\{y_i\}_{i\in I}$ be a sequence of integers, $\{a_{i}\}_{i \in I}$ a sequence of complex numbers and $\{b_i\}_{i \in I}$ a sequence of positive real numbers. When $|a_i| \leq b_i$ for each $i$ in $I$, the triangle inequality gives

\vspace{-3mm}
\begin{equation}
\int_{0}^{1} |\sum_{i \in I } a_i e(ty_i)|^2 \, dt \; = | \sum_{(i,j) \in I \times I} a_i {\bar a_j} \delta_{ij}\; | \; \leq \;  \sum_{(i,j) \in I \times I} b_ib_j \delta_{ij} \; = \; \int_{0}^{1} |\sum_{i \in I } b_i e(ty_i)|^2 \, dt\; ,
\end{equation}

\noindent
where $\delta_{ij} = 1$ when $y_i =y_j$ and 0 otherwise.

\end{nsc} 

\begin{nsc}{Counting Integer Points on Circles}

\noindent
When $a$ and $b$ are integers we write $(a,b) =1$ to mean that either one of $a$,$b$ is 1 and the other 0 or that $a$ and $b$ are both distinct from 0 and are coprime.

\begin{prop}
Let $H$ be a real number $\geq 1$, $m$ be a rational number and let $c_i$, $1 \leq i \leq 3$, be integers with $c_1 \neq 0$, $c_3 \geq 0$, $(c_1, c_2) =1$. Suppose that there are atleast three integer points $(x,y)$ satisfying $|x|, |y| \leq H$ and lying on the circle

\vspace{-3mm}
\begin{equation}
\label{c1}
(c_1 X - c_2)^2 + (c_1Y -m c_2)^2 = c_3 \; . 
\end{equation}

\noindent
If $m = \frac{p}{q}$ with $q > 0$ we then have

\vspace{-4mm}
\begin{equation}
\label{bound}
|c_1| \leq 4q(1+|m|)H, \hspace{2mm} |c_2| \leq 2qH^2 \hspace{1mm} \text{and} \hspace{1mm} c_3 \leq 36q^2(1+|m|)^2 H^4 \;. 
\end{equation}

\end{prop}

\vspace{4mm}
\noindent
{\sc Proof. ---} We follow the method of proof of Theorem 4, page 564 in  \cite{Hb}. Suppose that $p_i = (x_i,y_i)$, $1 \leq i \leq 3$, are three integer points on (\ref{c1}). Since the relation (\ref{c1}) is the same as 

\vspace{-4mm}
\begin{equation}
\label{c2}
{c_1}^{2} ({X}^2 + {Y}^2) -2c_1 c_2\, (X +mY) +(1 + m^2){c_2}^2 -  c_3 \; = 0\; , 
\end{equation}

\noindent
 we obtain the following relation of matrices on setting $(X,Y) = (x_i,y_i)$ in (\ref{c2}) for $1 \leq i \leq 3$.

\vspace{-2mm}
\begin{equation}
\label{mat}
\left(\begin{matrix}
{x_1}^2 + {y_1}^2 & x_1 +my_1 & 1 \\
{x_2}^2 + {y_2}^2 & x_2 +my_2 & 1 \\
{x_3}^2 + {y_3}^2 & x_3 +my_3 & 1 \\
\end{matrix}
\right)
\,
\left(\begin{matrix}
{c_1}^{2} \\
-2c_1 c_2 \\
(1+m^2){c_2}^2  - c_3
\end{matrix}\right) \; = \; 0 \, .
\end{equation}

\vspace{2mm}
\noindent
Let $M$ denote the $3 \times 3$ matrix and $c$ the vector $(c_1^2, -2c_1c_2, (1 + m^2)c_2^2 -c_3)$ in ${\bf Q}^{3}$ on the left hand side of (\ref{mat}). Since $c_1 \neq 0$, we have $c \neq 0$ and hence ${\rm det}(M) = 0$. Let $a = (a_1, a_2, a_3)$ be a solution distinct from 0 to  $Ma= 0$. Then the points $p_i$ lie on the conic

\vspace{-4mm}
\begin{equation}
\label{c3}
a_1(X^2 + Y^2) + a_2\, (X +mY) + a_3 \; = 0\; .
\end{equation}

\noindent
Since the points $p_i$ do not all lie on a line we must have $a_1 \neq 0$. Then (\ref{c3}) is an affine circle which intersects the affine circle (\ref{c1}) at the three points $p_i$. Since distinct affine circles intersect at no more than 2 points, we have $a = \alpha c$, for some complex number $\alpha$ which must necessarily be an element of ${\bf Q}$. Thus the rank of the matrix $M$ over ${\bf Q}$ is 2.

\vspace{2mm}
\noindent
Suppose now that the rows $i$ and $j$ of $M$, with $ i > j$, are linearly independent over ${\bf Q}$. Then $a = (a_1, a_2, a_3)$, where the $a_i$ are given by the relations

\vspace{-2mm}
\begin{equation}
\label{mat1}
a_1 = q\,{\rm det}\left(\begin{matrix}
x_i + my_i & 1 \\
x_j + my_j &1
\end{matrix}\right),
a_2 = -q\,{\rm det}\left(\begin{matrix}
{x_i}^{2} + y_i^{2} & 1 \\
x_j^{2} + y_j^{2} &1
\end{matrix}\right),
a_3 = q\,{\rm det}\left(\begin{matrix}
 x_i^{2} + y_i^{2} & x_i + my_i \\
 x_j^{2} + y_j^{2}  &x_j + my_j 
\end{matrix} \right) \; , 
\end{equation} 

\vspace{2mm}
\noindent
satisfies  $Ma = 0$. Since $a \neq 0$ and there is an $\alpha$ in ${\bf Q}$ such that $\alpha c = a$, we have 

\vspace{-3mm}
\begin{equation}
\frac{a_2}{a_1} = \frac{-2c_2}{c_1}, \hspace{3mm}  \frac{a_3}{a_1} = \frac{(1 + m^2)c_2^2 -c_3}{c_1^2} \;,
\end{equation}

\noindent
from which we deduce the following relations on setting $k = -2a_1/c_1 \neq 0$.

\vspace{-4mm}
\begin{equation}
\label{c2a}
k c_1 = -2a_1, \hspace{2mm} k c_2 = a_2, \; \; \text{and} \, \hspace{2mm} k^2 c_3 = (1+m^2)a_2^2 -4a_1 a_3 \;. 
\end{equation}

\noindent
Since the $a_i$ are integers and $(c_1, c_2) =1$, the relation $kc_2 = a_2$ shows that $c_1$ divides $-2a_1$ or that $k$ is an integer. Moreover, the triangle inequality applied to the relations in (\ref{mat1}) gives       

\vspace{-3mm}
\begin{equation}
\label{a2H}
|a_1|\, \leq \, 2q(1+|m|)H, \hspace{2mm} |a_2|\, \leq \, 2qH^2, \, \text{and} \hspace{2mm} |a_3| \leq 4q(1+|m|)H^3 \; .
\end{equation}

\noindent
We now obtain (\ref{bound}) on combining (\ref{c2a}) with (\ref{a2H}) using the triangle inequality and  $|k| \geq 1$.

\vspace{3mm}
\noindent
When $n$ is an integer $\geq 0$, $r(n)$ denotes the number of integer points $(x,y)$ satisfying $x^2 +y^2 = n$.

\vspace{2mm}
\begin{cor} Let $H$ be a real number $\geq  1$ and let $c_i$, $1 \leq i \leq 3$, be integers with $c_1 \neq 0$. The number of integer points $(x,y)$ satisfying $|x| \leq H$, $|y|\leq H$ and lying on
$(c_1X-c_2)^2 + (c_1Y-c_2)^2 = c_3$
does not exceed $\sup_{1 \leq n \leq 144H^4} r(n)$.

\end{cor}

\vspace{3mm}
\noindent
{\sc Proof. ---} Let $N(H)$ be the number of integer points satisfying the conditions of the corollary. We assume $N(H) \geq 1$ and set $d = c_1$ when $c_2 =0$ and $d = {\rm g.c.d.}(c_1, c_2)$ otherwise.  Then  $d^2$ divides $c_3$. Let $c_1/d = c_1^{*}$, $c_2/d = c_2^{*}$ and $c_3/d^2 = c_3^{*}$.  Then $N(H)$ is the same as the number of integer points $(x,y)$ satisfisying $|x| \leq H$, $|y| \leq H$ and lying on
$(c_1^{*}X-c_2^{*})^2 + (c_1^{*}Y-c_2^{*})^2 = c_3^{*} .$ 

\vspace{2mm}
\noindent
Plainly, $N(H) \leq r(c_3^{*})$. Since $c_1^{*} \neq 0$, $c_3^{*} \geq 0$ and $(c_1^{*}, c_2^{*}) =1$, Proposition 1 applied with $m=1$, $q=1$ implies that either $N(H) \leq 2$ or $c_3^{*} \leq 144H^4$ so that $N(H) \leq \sup(2, \sup_{1 \leq n \leq 144H^4} r(n))$, from which the corollary follows on noting that $\sup_{1 \leq n \leq 144H^4} r(n) \geq 2$ when $H \geq 1$.

\vspace{2mm}
\begin{cor} Let $I$ be a real interval of length $H \geq 1$. For all quadratic polynomials $P(T)$ in ${\bf Z}[T]$ and all integers $k$ the number of integer points $(x,y)$ in $I \times I$ lying on $P(X) + P(Y) = k$ does not exceed
$\sup_{1 \leq n \leq 144H^4} r(n)$.
\end{cor}

\vspace{3mm}
\noindent
{\sc Proof. ---} Suppose that $P(T) = a_0 T^2 + a_1 T + a_2$, with the $a_i$, $0 \leq i \leq 2$, integers and  $a_0 \neq 0$ and let $x_0$ be an integer in $I$. On completing the square and rearranging the terms we see that $P(x) + P(y) = k$ is equivalent to

\vspace{-4mm}
\begin{equation}
\label{nipc1}
(2a_0(x-x_0)+2a_0x_0 + a_1)^2 + (2a_0(y-x_0)+ 2a_0x_0 + a_1)^2 \; =\; 4a_0(k-2a_2)+2a_1^{2} \; ,
\end{equation}

\noindent
for any point $(x,y)$ in the real plane. Let us set $c_1 = 2a_0$, $c_2 = -(2a_0x_0+a_1)$ and write $c_3$ to denote the right hand side of (\ref{nipc1}). Since $x_0 \in I$, we have $|x-x_0| \leq H$ and $|y-x_0| \leq H$ for all $(x,y) \in I \times I$ and the number of integer points satisfying the conditions of the corollary does not exceed the number of integer points $(x,y)$ satisfying $|x| \leq H$, $|y| \leq H$ and lying on
$(c_1 X - c_2)^2 + (c_1Y -c_2)^2 = c_3$, so that the corollary follows from Corollary 1. 

\end{nsc} 

\begin{nsc}{A Variant of the Double Large  Sieve Inequality}

\noindent
The following lemma is the essence of the double large sieve inequality, modified slightly for our purpose. The proof follows pages 88 to 90, \cite{Graham} closely. 

\begin{le}
Let $\mathcal{X}$ be a $\delta$-spaced set of real numbers lying in the interval $[-P,P]$. Further, let $I$ be a finite set, $\{y_i\}_{i \in I}$ a sequence of integers and  $\{a_i\}_{i \in I}$ a sequence of complex numbers. When $T$ is a real number such that $|y_i| \leq T$ for all $i$ in $I$ we have the inequality 

\vspace{-3mm}
\begin{equation}
\label{ls1}
\left|\sum_{x \in \mathcal{X}}  f(x) \right|  \leq 
\pi \left( {\rm{Card}}(\mathcal{X})\, T + \frac{{\rm{Card}}(\mathcal{X})}{\delta} \right)^{1/2} 
(P + 2)^{1/2}  
\left(\int_{0}^{1} |f^{*}(t) |^2 \, dt \;\right)^{1/2} \, ,
\end{equation}

\noindent
where $f(t) = \sum_{i \in I } a_i e(ty_i)$ and $f^{*}(t) = \sum_{i \in I } |a_i| e(ty_i)$.

\end{le}

\vspace{3mm}
\noindent
{\sc Proof. ---} We have ${\rm Card}({\mathcal X}) \geq 2$ and that ${\mathcal X}$ is contained in $[-P,P]$. Therefore ${\rm Card}({\mathcal X}) -1 \geq {\rm Card}({\mathcal X})/2$ and , using (1) of 2.1, ${\rm Card}({\mathcal X}) -1 \leq 2P/\delta$. These relations show that (\ref{ls1}) holds when all the $y_i$ are 0. 

\vspace{2mm}
\noindent
Let us suppose that atleast one of the $y_i$ is distinct from 0. Since the $y_i$ are integers, we have $T \geq 1$. We then set $\epsilon = 1/4T$ and note that $\epsilon \leq 1$. Since $\phi_{\epsilon}(t) = \phi_{\epsilon}(-t)$, for any real number $x$ we have the relation 

\vspace{-3mm}
\begin{equation}
\label{ff}
\widehat{\phi_{\epsilon}}(y) e(xy) \; =\; \int_{\bf R} \phi_{\epsilon}(t-x) e(ty) \; dt.
\end{equation}

\noindent
From (\ref{ff}) we deduce that 

\vspace{-3mm}
\begin{equation}
\label{ff.5}
\sum_{x \in {\mathcal X}} f(x) \; = \; \sum_{x \in {\mathcal X}} \sum_{i \in I } a_i e(xy_i) \; = \; \int_{\bf R} \left(\sum_{x \in {\mathcal X}} \phi_{\epsilon}(t-x) \right) \left( \sum_{i \in I } \frac{a_i}{\widehat{\phi_{\epsilon}}(y_i)} e(ty_i) \right)\; dt.
\end{equation}

\noindent
Since ${\mathcal X}$ is a subset of $[-P,P]$ and $\epsilon \leq 1$, it follows that $\sum_{x \in {\mathcal X}} \phi_{\epsilon}(t-x)$ vanishes outside the interval $[-[P]-2, [P]+2]$. Let $\chi_{P}(t)$ denote the characteristic function of this interval. The Cauchy-Schwarz inequality then gives 

\vspace{-3mm}
\begin{equation}
\label{ff1}
\left| \sum_{x \in {\mathcal X}} f(x) \right| \; \leq \;  \left\|\sum_{x \in {\mathcal X}} \phi_{\epsilon}(t-x) \right\|_{2} \; \left\| \chi_{P}(t) \sum_{i \in I } \frac{a_i}{\widehat{\phi_{\epsilon}}(y_i)} e(ty_i) \right \|_{2} .
\end{equation}

\noindent
We have that $0 \leq \phi_{\epsilon} * \phi_{\epsilon}(t) \leq 2\epsilon$ for all $t$ in ${\bf R}$ and that the support of $\phi_{\epsilon} * \phi_{\epsilon}$ is $[-2\epsilon, 2\epsilon]$. These remarks together with $\phi_{\epsilon}(t) = \phi_{\epsilon}(-t)$  imply  

\vspace{-3mm}
\begin{equation}
\label{ff2.5}
\left\|\sum_{x \in {\mathcal X}} \phi_{\epsilon}(t-x) \right\|_{2}^{2} =
\sum_{(x, x^{\prime}) \in {\mathcal X} \times {\mathcal X}} \phi_{\epsilon}* \phi_{\epsilon}(x - x^{\prime}) \; \leq \; 2\epsilon \sum_{x \in \mathcal{X}} S(2\epsilon, x) \; \leq \;
2\epsilon\, {\rm{Card}}(\mathcal{X})\left( 1 + \frac{4\epsilon}{\delta}\right), 
\end{equation}

\noindent
where the last inequality follows from (1) of (2.1). Turning to the second term on the right hand side (\ref{ff1}), we note that since the $y_i$ are integers, $e(ty_i)$ is periodic of period 1 for each $i\in I$. Thus  

\vspace{-3mm}
\begin{equation}
\label{ff3}
\int_{-[P]-2}^{[P]+2} \left| \sum_{i \in I } \frac{a_i}{\widehat{\phi_{\epsilon}}(y_i)} e(ty_i) \right|^2 \; dt \, 
= 2([P] +2) \int_{0}^{1} \left| \sum_{i \in I } \frac{a_i}{\widehat{\phi_{\epsilon}}(y_i)} e(ty_i) \right|^2 \; dt \; 
%\leq \; 
%2([P] + 2) \left(\frac{\pi}{4\epsilon}\right)^{2} \int_{0}^{1} \left| \sum_{i \in I } a_i e(ty_i) \right|^2.
\end{equation}

\noindent
Recalling that $T = 1/4\epsilon$ and that the $y_i$ lie in the interval $[-T,T]$ we obtain $|a_i/\widehat{\phi_{\epsilon}}(y_i)| \leq \pi |a_i|/ 4 \epsilon$ for all $i \in I$ on using (\ref{phi1}) of (2.2). Using (\ref{ff3}) and the majorisation principle (2.3) we then conclude that 

\vspace{-3mm}
\begin{equation}
\label{ff4}
\left\| \chi_{P}(t)  \sum_{i \in I } \frac{a_i}{\widehat{\phi_{\epsilon}}(y_i)} e(ty_i) \right\|_{2}^2 \; dt \, \leq \,
2\left(\frac{\pi}{4\epsilon}\right)^{2} (P + 2) \int_{0}^{1} |f^{*}(t)|^2 \, dt.
\end{equation}

\noindent
The lemma now follows on combining (\ref{ff4}) with (\ref{ff1}) and (\ref{ff2.5}).

\vspace{3mm}
\begin{th}
Let $\mathcal{X}$ be a $\delta$-spaced set of real numbers lying in the interval $[-P,P]$. When $I$ is a finite set, $\mathcal{Y} = \{y_i\}_{i \in I}$ is a sequence of integers and  $\{a_i\}_{i \in I}$ are complex numbers we have the inequality

\vspace{-3mm}
\begin{equation}
\label{ls}
\sum_{x \in \mathcal{X}}  |\sum_{i \in I } a_i e(xy_i)|^2  \leq
\pi \left( {\rm{Card}}(\mathcal{X})\, \Delta(\mathcal{Y}) + \frac{{\rm{Card}}(\mathcal{X})}{\delta} \right)^{1/2}
(P + 2)^{1/2} \sup_{k} A_{\mathcal{Y}}^{1/2}(k) \| a\|^{2} \; ,
\end{equation}

\noindent
where $\Delta(\mathcal{Y})$ denotes $\sup_{(i,j) \in I \times I} |y_i -y_j|$, $A_{\mathcal{Y}}(k)$, for each integer $k$,  denotes the number of $(i,j) \in I \times I$ such that $y_i + y_j = k$ and $\| a\|^{2}$ denotes $\sum_{i \in I } |a_i|^2$.
  
\end{th}

\vspace{3mm}
\noindent
{\sc Proof. ---} Following the principle of Zhao \cite{zh}, we apply the preceding lemma  with $f(t)$ replaced by

\vspace{-4mm}
\begin{equation}
\label{idea}
g(t) \, =\, \sum_{(i,j) \in I \times I} a_i a_j e(t(y_i -y_j)) \; =\; |f(t)|^{2}
\end{equation}

\noindent
and with $T = \sup_{(i,j) \in I \times I} |y_i -y_j| = \Delta({\mathcal{Y}})$.  The theorem follows on noting that $g^{*}(t) = |f^{*}(t)|^2$ and using the estimate for $ \int_{0}^{1} |g^{*}(t)|^2 dt = \int_{0}^{1}|f^{*}(t)|^
4 dti \leq \sup_{k} A_{\mathcal{Y}}(k) \| a\|^{4}  
$.

\begin{cor}
Let $Q$ be a real number $\geq 1$ and suppose that ${\mathcal F}(Q)$ is the Farey sequence of order $Q$ and that $P(T) = c_0T^2 + c_1 T + c_2$ is a polynomial in ${\bf R}[T]$ with $c_0 \neq 0$, $c_1/c_0 = p/q \in {\bf Q} $ and $(p,q) =1$. When $\{a_i\}_{M < i \leq M+N}$ is a sequence of complex numbers indexed by the integers in the interval $(M, M+N]$, where $M$ and $N$ are integers with $N \geq 1$, we have

\begin{equation}
\label{ls0}
\sum_{x \in \mathcal{{\mathcal F}(Q)}}  \left| \sum_{M < i \leq M+N } a_i e( x P(i))\right|^2  \; \leq \; (Q^2 + Q \sqrt{c_0 N (|M| + 2N + |\frac{p}{q}| + 1)})\, \Pi \|a\|^2 \; ,
\end{equation}

\noindent
where \[
\Pi = \pi \left(\frac{2q}{c_0} + 1\right)^{\frac{1}{2}} \sup_{1 \leq n \leq 144 N^4} r(n).  
\]

\end{cor}

\noindent
Here $r(n)$ is the number of pairs of integers $(x,y)$ such that $x^2 + y^2 = n$.

\vspace{2mm}
\noindent
{\sc Proof. ---} We may assume $c_0 > 0$. We set $\alpha = \frac{c_0}{q}$ so that $P(T) = \alpha (qT^2 + pT) + c_2$. We take $I$ to be the set of integers in the interval $(M, M+N]$,

\vspace{-3mm}
\begin{equation}
\label{data}
y_i = qi^2 + pi, \; \; {\mathcal X} = \alpha {\mathcal F}(Q) \;\; 
\end{equation}

\noindent
We may then set $\delta = \frac{\alpha}{Q^2}$, ${\rm Card}({\mathcal X}) \leq Q^2$ and $P = \alpha$. We obtain Corollary 1 on applying Theorem 1 to the above data and taking into account that 

\vspace{-3mm}
\begin{equation}
\Delta(\mathcal{Y}) \leq |q|N(2N+|M|+1) + |p|N \leq qN( |M| + 2N + |\frac{p}{q}| +1)\; .
\end{equation}

\noindent
and that, by Corollary 2 to Proposition 1 of Section 3,  we have $A_{\mathcal{Y}}(k) \leq \sup_{1 \leq n \leq 144N^4} r(n)$, for all integers $k$.

\end{nsc}

\vspace{1cm}
\begin{flushleft}

{\em
Harish-Chandra Research Institute, \\
Chhatnag Road, Jhunsi,\\
Allahabad - 211 019, India.}\\

{\em email} : gyan@mri.ernet.in, suri@mri.ernet.in
\end{flushleft}

\end{document}